\documentclass[twoside]{article}

\usepackage{amssymb}
\usepackage{amsmath}
\usepackage{graphicx}

\makeatletter
\renewcommand{\@evenhead}{\thepage \hfil {\leftmark} \hfil}
\renewcommand{\@oddhead}{\hfil {\rightmark} \hfil \thepage}
\makeatother

\begin{document}
\author{\normalsize Hovhannes G. Tananyan, Rafayel R. Kamalian}
\title{\Large {\bf On Simultaneous 2-locally-balanced 2-partition for Two Forests with
Same Vertices}}
\date{} 

\maketitle
\begin{center}
{\small Yerevan State University}

{\small Russian-Armenian State University}

{\small e-mails: HTananyan@yahoo.com, rrkamalian@yahoo.com}
\end{center}

\begin{abstract}
The existence of a partition of the common set of the vertices
of two forests into two subsets, when difference of their capacities in the
neighbourhood of each vertex of each forest not greater than 2 is proved, 
and an example, which shows that improvement of the specified constant is
impossible is brought.
\end{abstract}

\bigskip

In this paper we continue researches started in [1-2], devoted to
locally-balanced partitions of a graph. We consider undirected graphs and
multigraphs without loops. The set of vertices of the multigraph $G$ is
denoted by $V(G)$, and the set of edges of $G$-by $E(G)$, and the maximum
degree of the vertices of $G$-by $\Delta(G)$. The eccentricity of a vertex
$v\in V(G)$ is denoted by $ex_{G}(v)$. Non-defined concepts can be found in
[3]. For $v\in V(G)$ we shall define sets $\gamma_{G}(v)=\{w\in V(G)/(w,v)\in
E(G)\}$ and $\eta_{G}(v)=\{e\in E(G)/v$ incident to $e\}$. A function
$f:M\rightarrow\{0,1\}$ is called $2$-partition of a finite set $M$. If $f$ is
a $2$-partition of a finite set $M$, then for $\forall M_{0}\subseteq M$ we
define the number $b_{f}(M_{0})$ as follows:

\begin{center}
$b_{f}(M_{0})=||\{m\in M_{0}/f(m)=1\}|-|\{m\in M_{0}/f(m)=0\}||$.
\end{center}

Let $G_{1}$and $G_{2}$ are undirected graphs without loops with $V(G_{1}%
)=V(G_{2})\equiv V$. The$\ 2$-partition $f$ of the set $V$ is called
simultaneous $k$-locally-balanced $(k\in Z,k\geq0)$ $2$-partition of the
graphs $G_{1}$and $G_{2}$ if:

\begin{center}
$\underset{i=1,2}{\max}$ $\underset{v\in V}{\max}b_{f}(\gamma_{G_{i}}(v))=k$.
\end{center}

Let $D$ is a tree, and let $v_{1}(D)\in V(D)$ is an arbitrarily chosen vertex.
For $i=0,1,...,ex_{D}(v_{1}(D))$ we define a subset $S_{i}\subseteq V(D)$ as follows:

\begin{center}
$S_{i}\equiv\{w\in V(G)/\rho(w,v_{1}(G))=i\}$.
\end{center}

For $i=1,2,...,ex_{D}(v_{1}(D))$ and $u\in S_{i-1}$ let's define
$S_{i}(u)\equiv\{w\in S_{i}/(w,u)\in E(D)\}$.

We define a family of subsets $X(D)$ of the set $V(D)$ as follows:

\begin{center}
$X(D)\equiv\{S_{i}(u)/1\leq i\leq ex_{D}(v_{1}(D)),u\in S_{i-1},S_{i}%
(u)\neq\varnothing\}\cup\{S_{0}\}$.
\end{center}

In the further we shall assume, that the consideration of any tree $D$ is
automatically implies the choice of the vertex $v_{1}(D)$.

Let $G$ is a forest, and $D_{1},D_{2},...,D\,_{k(G)}$ are its connected
components. Define a family of subsets $X(G)$ of the set $V(G)$ as follows:

\begin{center}
$X(G)\equiv\underset{i=1}{\overset{k(G)}{\cup}}X(D_{i})$.
\end{center}

Let $G_{1}$and $G_{2}$ are two forests with $V(G_{1})=V(G_{2})\equiv V$.
Define a bipartite multigraph $H(G_{1},G_{2})=(V_{1}(H(G_{1},G_{2}%
)),V_{2}(H(G_{1},G_{2})),E(H(G_{1},G_{2}))) $ as follows:

$V_{1}(H(G_{1},G_{2}))=X(G_{1})$,

$V_{2}(H(G_{1},G_{2}))=X(G_{2})$,

$E(H(G_{1},G_{2}))=\underset{v\in V}{\cup}\{(u,w)_{v}/u\in V_{1}(H(G_{1}%
,G_{2})),w\in V_{2}(H(G_{1},G_{2})),v\in u\cap w\}$,

where $E(H(G_{1},G_{2}))$ is understood as multiset containing different
elements like $(u,w)_{v_{1}}$ and $(u,w)_{v_{2}}$ with $v_{1}\neq v_{2}$ in a
case $|u\cap w|>1$.

It is not hard to see that for $\forall v\in V$

\begin{center}
$|\{(u,w)_{v}/u\in V_{1}(H(G_{1},G_{2})),w\in V_{2}(H(G_{1},G_{2})),v\in u\cap
w\}|=1$.
\end{center}

Taking into account that $G_{1}$ and $G_{2}$ are forests we can conclude from
the construction of the multigraph $H(G_{1},G_{2})$ that there exists an
one-to-one correspondence $\xi:V\rightarrow E(H(G_{1},G_{2}))$.

From the results of [4] it follows that there exists a $2$-partition $\varphi$
of the set $E(H(G_{1},G_{2}))$, at which for $\forall v\in V_{1}(H(G_{1}%
,G_{2}))\cup V_{2}(H(G_{1},G_{2}))$

\begin{center}
$b\varphi(\eta_{H(G_{1},G_{2})}(v))\leq1$.
\end{center}

\textbf{Theorem: }If $G_{1}$and $G_{2}$ are forests with $V(G_{1}%
)=V(G_{2})\equiv V$, then there exists a simultaneous $2$-locally-balanced
$2$-partition of $G_{1}$and $G_{2}$.

\textbf{Proof:} Define a $2$-partition $F$ of the set $V$ as follows: for
$\forall v\in V$ $\ F(v)\equiv\varphi(\xi(v))$. We shall be convinced that $F$
is a simultaneous $2$-locally-balanced $2$-partition of the forests $G_{1}$and
$G_{2}$. From the construction of the sets $X(G_{1})$ and $X(G_{2})$ it
follows that for $\forall v\in V$ \ $\exists A(v)\in X(G_{1})$ and $\exists
B(v)\in X(G_{2})$ such that $|\gamma_{G_{1}}(v)\backslash A(v)|\leq1$ and
$|\gamma_{G_{2}}(v)\backslash B(v)|\leq1$. Therefore it follows that
$b_{F}(\gamma_{G_{1}}(v))\leq$ $\ b_{F}(A(v))+1=$ $\ b_{\varphi}(\{e\in
E(H(G_{1},G_{2}))/\xi^{-1}(e)\in A(v)\})+1=$ $\ b_{\varphi}(\eta
_{H(G_{1},G_{2})}(v))+1\leq2$. Similarly, $b_{F}(\gamma_{G_{2}}(v))\leq2$.

Theorem is proved.

In the end we bring an example, which explains that not for arbitrary two
forests $G_{1}$and $G_{2}$ with $V(G_{1})=V(G_{2})\equiv V$ there exists a
$2$-partition $f$ of the set $V$, which is a simultaneous $k$-locally-balanced
$2$-partition of the forests $G_{1}$and $G_{2}$ for $k\leq1$.

\textbf{Example:} Define trees $G_{1}$and $G_{2}$ as follows:

$G_{1}=(\{v_{1},v_{2},v_{3},v_{4},v_{5}\},\{(v_{1},v_{2}),(v_{2},v_{3}%
),(v_{3},v_{4}),(v_{4},v_{5})\})$,

$G_{2}=(\{v_{1},v_{2},v_{3},v_{4},v_{5}\},\{(v_{1},v_{5}),(v_{5},v_{4}%
),(v_{3},v_{2}),(v_{2},v_{1})\})$.

\begin{center}
\begin{figure}
\begin{center}
\includegraphics[natheight=1.260000in, natwidth=2.874600in, height=1.2955in, width=2.9187in]{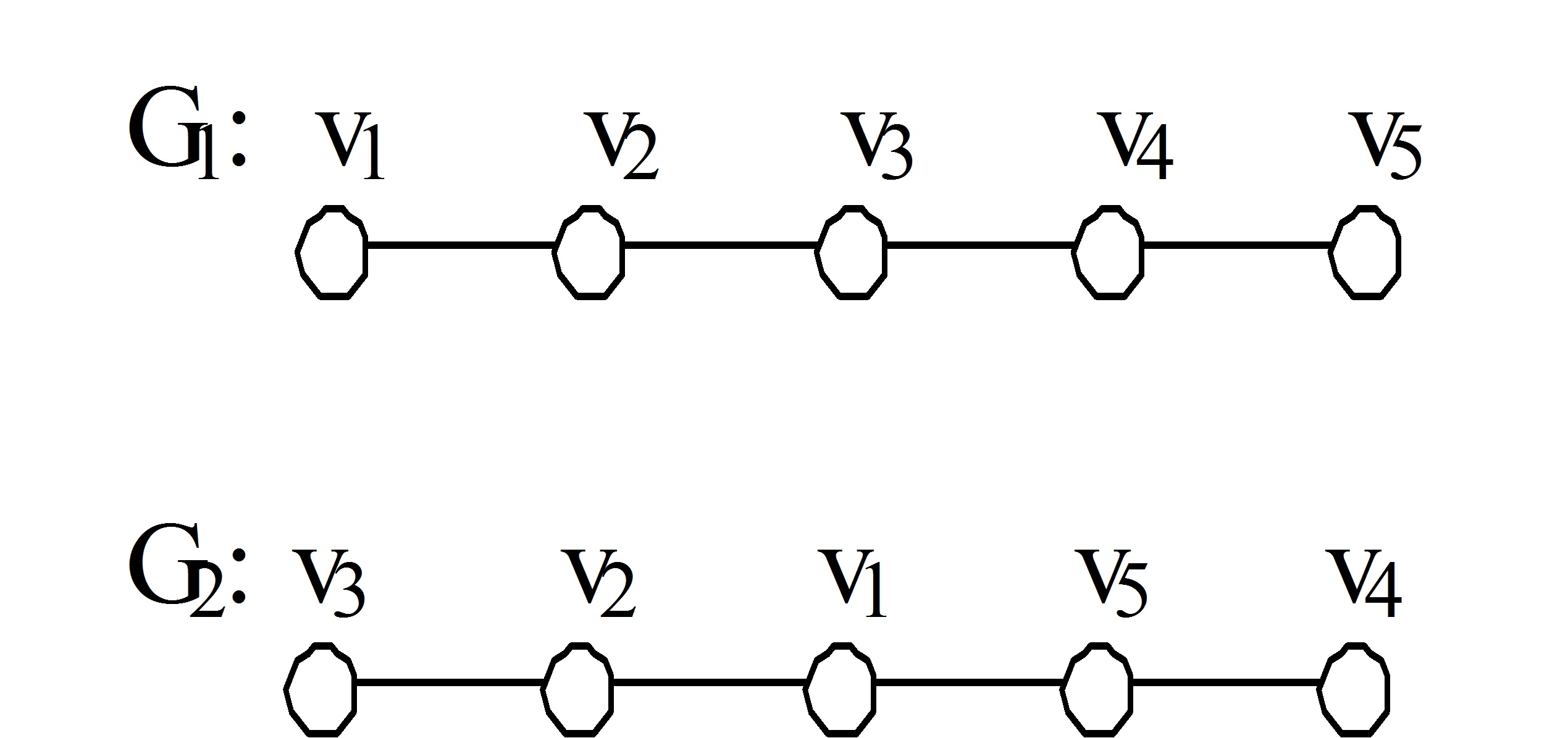}
\end{center}
\end{figure}
\end{center}

Let's assume that there exists a $2$-partition $f$ of the set $\{v_{1}%
,v_{2},v_{3},v_{4},v_{5}\}$, which is a simultaneous $k$-locally-balanced
$2$-partition of the trees $G_{1}$and $G_{2}$ for $k\leq1$. Without
restriction of a generality we can suppose, that $f(v_{1})=0$. From
$\gamma_{G_{1}}(v_{2})=\{v_{1},v_{3}\}$ and $\gamma_{G_{1}}(v_{4}%
)=\{v_{3},v_{5}\}$ we can conclude that $f(v_{3})=1$ and $f(v_{5})=0$ . Hence,
from $\gamma_{G_{2}}(v_{5})=\{v_{1},v_{4}\}$ and $\gamma_{G_{2}}%
(v_{1})=\{v_{2},v_{5}\}$ we can conclude that $f(v_{4})=1$ and $f(v_{2})=1$.
But it means that $b_{f}(\gamma_{G_{1}}(v_{3}))=2$, which contradicts the
property of $f$.

\end{document}